\documentclass{amsart}

\usepackage{amsmath}
\usepackage{amsfonts}
\usepackage{graphics}
\usepackage{epsfig} 
\usepackage{amssymb}
\usepackage{amscd}
\usepackage[all]{xy}
\usepackage{latexsym}
\usepackage{graphicx}
\usepackage{multirow}
\usepackage{geometry}

\newtheorem{thm}{Theorem}[section]

\newtheorem{prop}[thm]{Proposition}

\theoremstyle{definition}

\newtheorem{rem}[thm]{Remark}
\newtheorem*{ack}{Acknowledgement}

\numberwithin{equation}{section}
\numberwithin{figure}{section}

\def\Hom{{\text{\rm{Hom}}}}

\def\rchi{{\hbox{\raise1.5pt\hbox{$\chi$}}}}
\def\Aut{{\text{\rm{Aut}}}}

\def\a{\alpha}
\def\b{\beta}
\def\lam{\lambda}
\def\Lam{\Lambda}
\def\gam{\gamma}
\def\Gam{\Gamma}

\newcommand{\bP}{{\mathbb{P}}}
\newcommand{\bC}{{\mathbb{C}}}

\newcommand{\bN}{{\mathbb{N}}}
\newcommand{\bE}{{\mathbb{E}}}

\newcommand{\bQ}{{\mathbb{Q}}}

\newcommand{\cM}{{\mathcal{M}}}

\newcommand{\cH}{{\mathcal{H}}}

\newcommand{\la}{{\langle}}
\newcommand{\ra}{{\rangle}}
\newcommand{\half}{{\frac{1}{2}}}

\newcommand{\hxi}{{\hat{\xi}}}
\newcommand{\hatH}{{\widehat{\cH}}}

\begin{document}
\large
\setcounter{section}{0}

\title[Polynomial recursion for 
{Hodge} integrals]{Polynomial recursion formula for 
 linear {Hodge} integrals}
\author[M.~Mulase]{Motohico Mulase}  
\address{
Department of Mathematics\\
University of California\\
Davis, CA 95616--8633}
\email{mulase@math.ucdavis.edu}
\author[N.~Zhang]{Naizhen Zhang}  
\address{
Department of Mathematics\\
University of California\\
Davis, CA 95616--8633}
\email{nzhzhang@math.ucdavis.edu}

\begin{abstract}
We establish a polynomial recursion formula for
linear Hodge integrals. It is obtained as 
the Laplace 
transform of the cut-and-join equation for the simple 
Hurwitz numbers. We 
show that the recursion recovers the Witten-Kontsevich
theorem when restricted to the top degree terms, and 
also the combinatorial factor of the 
$\lam_g$ formula as the lowest degree terms.
\end{abstract}

\subjclass[2000]{14H10, 14N10, 14N35; 05A15, 05A17; 81T45}

\maketitle

\allowdisplaybreaks

\textit{\centerline{Dedicated to Herbert Kurke on
the occasion of his
70th birthday}}

\tableofcontents

\section{Introduction}
\label{sect:intro}

The purpose of this paper is to establish a
topological recursion formula for linear Hodge integrals
in terms of \emph{polynomial} generating functions.
Let $\overline{\cM}_{g,\ell}$ be the Deligne-Mumford
moduli stack of stable curves of genus $g$ and $\ell$ distinct
marked points  subject to $2g-2+\ell>0$.
We denote by $\psi_i$  the $i$-th cotangent class
of $\overline{\cM}_{g,\ell}$,
and by $\lam_j=c_j(\bE)$ the $j$-th Chern class of 
the Hodge bundle $\bE$ on $\overline{\cM}_{g,\ell}$.
By \emph{linear} Hodge integrals we mean the
rational numbers
$$
\la \tau_{n_1}\cdots\tau_{n_\ell}\lam_j\ra_{g,\ell}
=
\int_{\overline{\mathcal{M}}_{g,\ell}}
\psi_1 ^{n_1}\cdots \psi_\ell ^{n_\ell}\lam_j.
$$
Following \cite{EMS, GJV2} we define a series of polynomials by 
a recursion formula
$$
\hxi_{n+1}(t) = t^2(t-1)\frac{d}{dt}\hxi_n(t)=D\hxi_n(t)
$$
with the initial condition $\hxi_0(t) = t-1$.
The differential operator $D=t^2(t-1)\frac{d}{dt}$
found in \cite[Example~4.1]{GJV2} 
simplifies many of  the combinatorial difficulties of the
linear Hodge integrals and Hurwitz numbers.
The degree of $\hxi_n(t)$ is $2n+1$. 
We consider symmetric polynomials of degree
$3(2g-2+\ell)$,
\begin{equation}
\label{eq:hatH}
\hatH_{g,\ell}(t_1,\dots,t_\ell)
=\sum_{n_1,\dots,n_\ell}
\la \tau_{n_1}\cdots\tau_{n_\ell}\Lam_g ^\vee(1)\ra_{g,\ell}
\prod_{i=1} ^\ell \hxi_{n_i}(t_i),
\end{equation}
where $\Lam_g ^\vee(1)
=1-\lam_1+\cdots+(-1)^g \lam_g$. 
The following is our main theorem.

\begin{thm}
\label{thm:main}
The polynomial generating functions of the
linear Hodge integrals {\rm{(\ref{eq:hatH})}} satisfy the
following topological recursion formula
\begin{multline}
\label{eq:main}
\left(
2g-2+\ell +\sum_{i=1} ^\ell
\frac{1}{t_i}
D_i
\right)
\hatH_{g,\ell}(t_L)
\\
=
\sum_{i< j}
 \frac{  t_i ^2\hxi_0(t_j)D_i
    \hatH_{g,\ell-1}\left(t_{L\setminus \{j\}}\right)
    -
      t_j ^2\hxi_0(t_i)D_j
    \hatH_{g,\ell-1}\left(t_{L\setminus \{i\}}\right)}{t_i-t_j}
\\
+
\sum_{i=1} ^\ell
\left[
D_{u_1}D_{u_2}
\hatH_{g-1,\ell+1}\left(u_1,u_2,t_{L\setminus\{i\}}\right)
\right]_{u_1=u_2=t_i}
\\
+
\half
\sum_{i=1} ^\ell
\sum_{\substack{g_1+g_2 = g\\
J\sqcup K= L\setminus\{i\}}} ^{\rm{stable}}
D_i
\hatH_{g_1,|J|+1}(t_i,t_J)\cdot 
D_i
\hatH_{g_2,|K|+1}(t_i,t_K) ,
\end{multline}
where $D_i = t_i ^2 (t_i-1)\frac{\partial }{\partial t_i}$.
The last summation is taken over all partitions
$g=g_1+g_2$  of the genus $g$ and disjoint union  decompositions
$J\sqcup K= L\setminus\{i\}$
satisfying the stability conditions
$2g_1-1+|J|>0$ and $2g_2-1+|K|>0$. Here 
$L=\{1,2,\dots,\ell\}$ is the index set, and for 
a subset $I\subset L$ we write $t_I = (t_i)_{i\in I}$. 
\end{thm}

The recursion formula
(\ref{eq:main}) is a \emph{topological} recursion in the sense that
it gives the generating function 
of linear Hodge integrals of complexity $2g-2+\ell=n$ in terms of those
of complexity $n-1$. The same topological structure appears
in other recursion formulas
such as those discussed in 
\cite{Dijkgraaf, DVV, EO1, EO2, GJV2, LX, LX2, Mir1, Mir2, MS}.

We prove Theorem~\ref{thm:main} by computing the
Laplace transform of the Hurwitz number $h_{g,\mu}$
as a \emph{function}
of a partition $\mu$. 
Let $f:X\rightarrow \bP^1$ be a morphism of 
connected nonsingular algebraic curve $X$ of genus $g$
onto the projective line defined over $\bC$. If we regard 
$f$ as a meromorphic function on $X$, then the \emph{profile}
of $f$ is the list of orders of its poles being considered as
a \emph{partition} of the degree of $f$. The \emph{Hurwitz}
number $h_{g,\mu}$  we deal with in this paper is the number of 
topological types of $f$ of given genus $g$ and profile
$\mu$ being counted with the weight $1\big/|\Aut(f)|$. 
The celebrated \emph{cut-and-join 
equation} of Goulden, Jackson, and Vakil \cite{GJ, V}
(which was essentially known to Hurwitz \cite{H})
applied to the Laplace transformed Hurwitz numbers 
is exactly the polynomial recursion (\ref{eq:main}).
The idea of taking the Laplace transform of the
cut-and-join equation comes from \cite{EMS}.
It is shown in \cite{EMS} that (\ref{eq:main}) implies
the Bouchard-Mari\~no conjecture on the topological 
recursion for Hurwitz numbers
\cite{BM}, which is the simplest case of the 
more general conjecture on the closed and open
Gromov-Witten invariants of toric Calabi-Yau 3-folds
 \cite{BKMP}.

The significance of (\ref{eq:main}) being
a polynomial is two-fold. 
Firstly, the leading coefficients of $\hatH_{g,\ell}$ are
the $\psi$-class intersection numbers. 
It was proved by Okounkov and Pandharipande 
\cite{OP1} that the large partition asymptotics
of the Hurwitz numbers recover the Witten-Kontsevich
theorem, i.e., the Virasoro
constraint condition for the $\psi$-class intersection 
numbers \cite{DVV, K1992, W1991}. 
Since the Laplace transform contains more information 
than the asymptotic
behavior, the proof of the Witten conjecture \cite{W1991}
becomes just comparing the leading 
coefficients of the polynomial equation (\ref{eq:main}).
The second significance is that the coefficients
of the \emph{lowest} degree
terms are the linear Hodge integrals containing 
the $\lam_g$-class. The topological recursion 
recovers the formula
for 
$
\la \tau_{n_1}\cdots\tau_{n_\ell}\lam_g\ra_{g,\ell}
$
in terms of 
$
\la \tau_{2g-1}\lam_g\ra_{g,1}
$.
We remark that the same \emph{polynomiality} 
is observed in \cite{Kazarian, KazarianLando}
in the context of integrable systems.

We note that all the formulas in this paper have been  
more or less established in various different formulations
\cite{CLL, GJV1, GJV2, GJV3,
Kazarian, LZZ}.
Since (\ref{eq:main}) is \emph{equivalent} to the cut-and-join
equation, logically speaking one may say there is nothing new. The contribution
of this paper is
the simple expression of our formulation of the
cut-and-join equation (\ref{eq:caj})
and a new point of view of understanding (\ref{eq:main})
as the Laplace transform of (\ref{eq:caj}). 
It gives a clear and unified picture of some of the 
results established in \cite{CLL, GJV3, Kazarian}.

The paper is organized as follows. We begin with
setting our notations and 
reviewing definitions of 
Hurwitz numbers in Section~\ref{sect:Hurwitz}. 
In Section~\ref{sect:CAJ} we formulate  the cut-and-join
equation as a functional equation for functions
in partitions. Although there are a large number of
literature on the subject \cite{CLL, GJ, GJV1, GJV2, GJV3,
Kazarian, KazarianLando, LZZ, LLZ, V, Zhou}, we provide a 
full detail in this section because we wish to arrive at a 
simpler formulation of the equation.
 We then introduce the idea
of Laplace transformation following \cite{EMS}
in Section~\ref{sect:LT}. Here the role of the
Lambert curve, the \emph{spectral curve} 
of the topological recursion for Hurwitz numbers
introduced in \cite{BEMS, BM, EMS, EO1},
is identified as the \emph{Riemann surface}
of  a meromorphic function that is obtained by
the Laplace transform. The following Section~\ref{sect:recursion}
establishes Theorem~\ref{thm:main}. 
In the final section we derive the Dijkgraaf-Verlinde-Verlinde
formula \cite{DVV} for the Witten-Kontsevich
theorem \cite{K1992, W1991} from (\ref{eq:main})
as a simple corollary. We also give the combinatorial 
coefficient of the $\lam_g$ formula
\cite{FP1,FP2} from the topological recursion.

\begin{ack}
The authors thank the American Institute of 
Mathematics for the hospitality during their stay that
promoted this collaboration. They are 
grateful to Ravi Vakil, Lin Chen, and the referee for useful 
comments. 
M.M.\ thanks Herbert Kurke for giving him the
opportunity to lecture on Hurwitz numbers 
based on \cite{OP1, OP2} at
Humboldt Universit\"at zu Berlin in 2002 and 2005.
M.M.\ also thanks the NSF, Kyoto 
University,
 the Institute for the Physics 
and Mathematics of the Universe in Tokyo,
the Osaka City University Advanced Mathematical 
Institute,    T\^ohoku University, KIAS in Seoul,
and the University of Salamanca 
for their hospitality and financial support during the 
preparation 
of this work.
\end{ack}

\section{Hurwitz numbers}
\label{sect:Hurwitz}

Let $X$ be a nonsingular complete algebraic curve of genus $g$
defined over the complex number field $\bC$,
and 
 $f:X\rightarrow \bP^1$ a morphism 
 of $X$ to the projective line
  $\bP^1$. If we regard $f$ a meromorphic
function on the Riemann surface
$X$, then the inverse image $f^{-1}(\infty) = \{p_1,\dots,
p_\ell\}$ of
$\infty\in\bP^1$ is the set of poles of $f$. We can name
these $\ell$ points   so that the list of pole orders  becomes a \emph{partition}
$\mu = (\mu_1\ge \mu_2 \ge
\cdots \ge \mu_\ell >0)$ of  the degree of the map. 
Thus the \emph{size} of this partition
$
|\mu| =
\mu_1 + \cdots + \mu_\ell
$
is
 $\deg f $, and its \emph{length}  $\ell(\mu) = \ell$
is the number of poles of $f$. 
Each part $\mu_i$ determines a local description of the map $f$,
which is given by $z\longmapsto z^{\mu_i}$ in terms of a local
coordinate $z$ of $X$ around $p_i$. 
A critical point, or a \emph{ramification point},
 of $f$ is a point $p\in X$ at which 
the derivative vanishes $df(p)=0$, and $w=f(p)$ is a critical 
value, or
a \emph{branched point} of $f$. 
Let $B\subset \bP^1$ be the set of all branched points  of $f$.
Then 
\begin{equation}
\label{eq:covering}
f\big|_{f^{-1}(\bP^1\setminus B)}:
f^{-1}(\bP^1\setminus B)\longrightarrow
\bP^1\setminus B
\end{equation}
is a topological covering of degree $|\mu|$. 
When the derivative $df$ has a {simple} zero
at $p$, we say $p$ is a \emph{simple ramification point}
of $f$. If over every branched point except for
$\infty$ there is exactly one
simple ramification point, then we call $f$ 
a \emph{Hurwitz cover}. The partition $\mu$ gives the \emph{profile} of
a Hurwitz cover. 
The number 
$h_{g,\mu}$ of topological
types of Hurwitz covers of given genus $g$ and profile 
$\mu$, counted with the weight factor $1/|\Aut f|$, is the \emph{Hurwitz
number} we are interested in this paper. 
To be more precise, we study
 $h_{g,\mu}$ as a \emph{function}
of partition $\mu$. We will compute
the \emph{Laplace transform} of $h_{g,\mu}$ and
find the equations that they satisfy.

Let $r$ denote the number of simple ramification 
points of $f$. This gives the dimension of the Hurwitz scheme,
i.e., the moduli space of all Hurwitz covers for a given
genus and a profile
\cite[Section~7.3.2]{OP1}. Since (\ref{eq:covering}) is a topological 
covering, the Euler characteristic of $f^{-1}(\bP^1\setminus B)$
is given by
$$
\rchi\big(f^{-1}(\bP^1\setminus B)\big) = \deg f \cdot
\rchi(\bP^1\setminus B) = |\mu| (1-r).
$$
On the other hand, since $f^{-1}(x)$ contains 
exactly $\deg f -1$ points for every $x\in B\setminus\{\infty\}$
and since $f^{-1}(\infty)$ has $\ell$ points, 
$$
\rchi\big(f^{-1}(\bP^1\setminus B)\big)=
2-2g(X) -\ell - r(|\mu|-1) .
$$
We thus obtain the \emph{Riemann-Hurwitz formula}
\begin{equation}
\label{eq:RH}
r = r(g,\mu) = 2g - 2 +\ell +|\mu|.
\end{equation}
The celebrated Ekedahl-Lando-Shapiro-Vainshtein
 formula \cite{ELSV, GV1, OP1} relates 
Hurwitz numbers and linear Hodge integrals 
on the Deligne-Mumford moduli stack 
$\overline{\mathcal{M}}_{g,\ell}$
consisting of stable algebraic curves of genus $g$ with $\ell$ distinct
nonsingular marked points subject to the
stability condition $2g-2+\ell >0$. 
Denote by $\pi_{g,\ell}:\overline{\mathcal{M}}_{g,\ell+1}
\rightarrow \overline{\mathcal{M}}_{g,\ell}$
the natural projection and by $\omega_{\pi_{g,\ell}}$ the relative
dualizing sheaf of the universal curve $\pi_{g,\ell}$. 
The \emph{Hodge} bundle $\bE$ 
on $\overline{\mathcal{M}}_{g,\ell}$
is defined by $\bE = (\pi_{g,\ell})_* \omega_{\pi_{g,\ell}}$,
and  the $\lambda$-classes 
are  the Chern classes 
$$
\lambda_i = c_i(\bE)
\in H^{2i}(\overline{\mathcal{M}}_{g,\ell}, \bQ)
$$
of the Hodge bundle.
Let $\sigma_i:\overline{\mathcal{M}}_{g,\ell}
\rightarrow \overline{\mathcal{M}}_{g,\ell+1}$
be the $i$-th tautological section of $\pi$, and put
$\mathcal{L}_i = \sigma_i ^*(\omega_{\pi_{g,\ell}})$. 
The $\psi$-classes are defined by
$$
\psi_i = c_1(\mathcal{L}_i) \in 
H^{2}(\overline{\mathcal{M}}_{g,\ell}, \bQ).
$$
The \emph{linear Hodge integrals}
are rational numbers defined by
$$
\la \tau_{n_1}\cdots\tau_{n_\ell}\lam_j\ra_{g,\ell}
=
\int_{\overline{\mathcal{M}}_{g,\ell}}
\psi_1 ^{n_1}\cdots \psi_\ell ^{n_\ell}\lam_j,
$$
which are $0$ unless $n_1+\cdots+n_\ell +j = 3g-3+\ell$. 
Let us denote by
${\Lambda_g^{\vee}(1)}= 1-\lam_1+\cdots+(-1)^g \lam_g$. 
The ELSV formula states
\begin{equation}
\label{eq:ELSV}
	h_{g,\mu} = \frac{r(g,\mu)!}{|\Aut( \mu)|}\;  
	\prod_{i=1}^{\ell(\mu)}\frac{\mu_i ^{\mu_i}}{\mu_i!}
	\int_{\overline{\mathcal{M}}_{g,\ell(\mu)}} 
	\frac{\Lambda_g^{\vee}(1)}
	{\prod_{i=1}^{\ell(\mu)}\big( 1-\mu_i \psi_i\big)},
\end{equation}
where  $\Aut(\mu)$ is 
the permutation group that interchanges the equal parts of 
$\mu$. 
The appearance of this automorphism factor is due to the
difference between giving a profile $\mu$ and 
\emph{naming} all points in $f^{-1}(\infty)$. 
If all parts of $\mu$ are distinct, then the poles of 
$f$ are naturally labeled by the pole order. But when
two or more parts are the same, there is no way to
distinguish the Hurwitz covers obtained by
interchanging these
poles of the same order. The factor $1\big/|\Aut(\mu)|$
takes care of this overount.

Although $\overline{\mathcal{M}}_{g,\ell}$
is defined as the moduli stack of \emph{stable} curves satisfying the
stability condition
$2-2g-\ell <0$, Hurwitz numbers 
are well defined for \emph{unstable} geometries
 $(g,\ell) = (0,1)$ and $(0,2)$. 
It is an elementary exercise to 
show that
\begin{equation*}
h_{0,k} = k^{k-3}\qquad \text{and}\qquad
h_{0,(\mu_1,\mu_2)} = \frac{(\mu_1+\mu_2)!}{\mu_1+\mu_2}\cdot
\frac{\mu_1^{\mu_1}}{\mu_1!}\cdot \frac{\mu_2^{\mu_2}}{\mu_2!}.
\end{equation*}
The ELSV formula remains true for unstable cases
by \emph{defining}
\begin{align}
\label{eq:01Hodge}
&\int_{\overline{\cM}_{0,1}} \frac{\Lambda_0 ^\vee (1)}{1-k\psi}
=\frac{1}{k^2},\\
\label{eq:02Hodge}
&\int_{\overline{\cM}_{0,2}} 
\frac{\Lambda_0 ^\vee (1)}{(1-\mu_1\psi_1)(1-\mu_2\psi_2)}
=\frac{1}{\mu_1+\mu_2}.
\end{align}

\section{The cut-and-join equation}
\label{sect:CAJ}

The Hurwitz numbers satisfy a set of combinatorial equations
called the \emph{cut-and-join} equation 
of 
\cite{GJ, V}. It is essentially the same relation Hurwitz
dealt with in his seminal paper \cite{H}. Due to the 
modern formulation in these more recent papers,
 the combinatorial equation  has become an effective
tool of algebraic geometry for studying  Hurwitz 
numbers and many related subjects \cite{CLL, GJV1, GJV2, GJV3, GV1,GV2,KazarianLando,
 LLZ, OP1, Zhou}. 
In  this section we review the equation
following \cite{GJ, LZZ,V, Zhou}, and give its simplest
formulation that is suitable to compute its Laplace
transform in Section~\ref{sect:recursion}.

 The topological covering (\ref{eq:covering})
gives rise  to a unique point
 in the \emph{character variety}
 \begin{equation}
 \label{eq:char}
 \rho\in \Hom\big(\pi_1(\bP^1\setminus B),S_d\big)\big/S_d,
 \end{equation}
where $S_d$ is the symmetric group of $d=|\mu|$ letters
and its action on the set of homomorphisms is through 
conjugation. Since the character variety classifies 
\emph{all} topological coverings, we need to determine
the condition for a covering to be a Hurwitz cover. 
Let us list the $r+1$ points
in $B$  as
$$
B = \{x_1, \dots, x_r, \infty\}.
$$
Choose a base point $*$ on $\bP^1\setminus B$, and
denote by $\gam_k$ a closed path starting from $*$ that goes
around $x_k$ in the positive direction, and comes back 
to $*$. The loop $\gam_\infty$ is the loop going around
$\infty$. Then up to conjugation, we have
$$
\pi_1(\bP^1\setminus B) \cong
\la \gam_1,\dots,\gam_r,\gam_\infty\;|\;
\gam_1\cdots\gam_r\cdot \gam_\infty = 1\ra.
$$
Now recall that over each $x_k$ there is only one
ramification point, say $p_k$, which is simple. 
Therefore, in terms of the  representation $\rho$ corresponding
to the Hurwitz cover $f$, the generator $\gam_k$ is mapped to a 
transposition $(ab)\in S_d$. 
Next, recall that the ramification behavior over $\infty$ is
determined by the profile $\mu$, and that each part $\mu_i$
determines the map $f$ locally as $z\longmapsto z^{\mu_i}$.
In terms of the representation, this means that 
$$
\rho(\gam_\infty) = c_1c_2\cdots c_\ell,
$$
where 
$$
c_1\sqcup\cdots\sqcup c_\ell = \{1,2,\dots,d\}
$$
is a disjoint cycle decomposition of the index set 
 and each $c_i$ is a cycle of 
length $\mu_i$.

The cut-and-join equation represents the number of
Hurwitz covers of a given genus $g$ and profile $\mu$
in terms of those with profiles obtained by either 
\emph{cutting} a part into two pieces, or \emph{joining}
two parts together. 
Let $p\in X$ be a point at which the covering $f:X\rightarrow \bP^1$
is simply ramified. Locally we can name sheets, so we assume
sheets $a$ and $b$ are ramified over $x_r=f(p)\in B\subset \bP^1$. 
In terms of the representation we have 
$\rho(\gam_r) = (ab)\in S_d$. 
When we merge $x_r$ to $\infty$, the generators
$\gam_r$ and $\gam_\infty$ of 
$\pi_1(\bP^1\setminus B)$ are replaced by their product 
$\gam_r\gam_\infty$. The representation $\rho$ maps 
this generator to $(ab)c_1\cdots c_\ell$. 
Now one of the two things happen:

\begin{enumerate}

\item The \emph{cut case}, in which both sheets are ramified at
  the same point
$p_i$
of the inverse image $f^{-1}(\infty) = \{p_1,\dots, p_\ell\}$.
In terms of $\rho$, this means both indices $a$ and $b$ 
are contained in the same cycle $c_i$. 
Since $c_1,\dots, c_\ell$ are disjoint, we only need to 
calculate $(ab)c_i$. By re-naming 
all the sheets and assuming $a<b=a+\a<\mu_i=\a+\b$, we can compute
\begin{align*}
&\big(a[a+\a]\big)\big(12\cdots [a-1] a [a+1]\cdots [a+\a]\cdots 
[\a+\b]\big)
\\
&=
\big(a[a+1]\cdots [a+\a-1]\big)
\big([a+\a][a+\a+1]\cdots [\a+\b]12\cdots[a-1]\big).
\end{align*}
 The result is the product of two disjoint cycles of length
 $\a$ and $\b$. 
Thus the merging eliminates
a profile $\mu$ and creates a new profile
$$
(\mu_1,\dots,\widehat{\mu_i},\dots,  \mu_{\ell},\alpha,\b)
=\big(\mu(\hat{i}),\a,\b\big)
$$
 of length $\ell +1$. 
 Here the $\widehat{\;\;}$ sign means removing the entry. 
 Note that the size of the partition $|\mu|$
 is unchanged,
 because it is the degree of the map $f$. When $\a$ is chosen, 
 the  total number of 
 such cuttings is $\a+\b$ because this is the number of 
 choices for $a$ in the index set $\{1,2,\dots,\a+\b\}$. 
 We also note that when $\a=\b$, the number is actually $\a$,
 instead of $\a+\b$.

\item   The
\emph{join case}, in which
 sheets $a$ and $b$ are ramified at two distinct
points, say $p_i$ and $p_j$, above $\infty$. In other words,
$a\in c_i$ and $b\in c_j$. Again by re-numbering,
we can calculate
\begin{align*}
&(ab)\big(12\cdots [a-1] a [a+1]\cdots \mu_i\big)
\big([\mu_i+1]\cdots[b-1] b [b+1]\cdots [\mu_i+\mu_j]\big)
\\
&=\big(12\cdots[a-1]b[b+1]\cdots [\mu_i+\mu_j][\mu_i+1]
\cdots[b-1]a[a+1]\cdots\mu_i\big).
\end{align*}
Thus the result of merging 
creates a new profile
$$
(\mu_1,\dots,\widehat{\mu_i},\dots,\widehat{\mu_j},\dots,\mu_\ell,
\mu_i+\mu_j) 
=\big(\mu(\hat{i},\hat{j}) ,\mu_i+\mu_j\big)
$$
of length $\ell-1$ and size
$|\mu|$. The total number of ways to 
make the join is $\mu_i\mu_j$, because we have
$\mu_i$-choices for $a$ and $\mu_j$-choices for $b$. 
\end{enumerate}

To utilize the above consideration into 
Hurwitz numbers, let us introduce
the generating function 
of Hurwitz numbers 
\begin{equation}
\label{eq:anotherHurwitz}
\mathbf{H}(s,\mathbf{p})=\sum_{g\ge 0}\sum_{\ell\ge 1} \mathbf{H}_{g,\ell}
(s,\mathbf{p});
\qquad
\mathbf{H}_{g,\ell}
(s,\mathbf{p}) =\sum_{\mu: \ell(\mu)=\ell} h_{g,\mu}\mathbf{p}_{\mu} \frac{s^{r(g,\mu)}}{r(g,\mu)!},
\end{equation}
where $\mathbf{p}_{\mu} = p_{\mu_1}p_{\mu_2}\cdots p_{\mu_{\ell}}$,
and $r(g,\mu)$ is the number of simple ramification points
(\ref{eq:RH}).
The summation in $\mathbf{H}_{g,\ell}
(s,\mathbf{p})$  is over all partitions of length $\ell$. 
Here $p_1,p_2,p_3,\dots$ are parameters
that encode the information of partitions. 
The other parameter $s$ counts the number $r$ of simple 
ramification points. Since $r$ and $\mu$ recover
the genus $g$, $s$ is a \emph{topological} parameter. 
Note that 
merging $x_r$ to $\infty$ means decreasing $r$ by $1$, or 
differentiating the generating function with respect to 
$s$. The result of this differentiation is
the cut and join operations discussed above. 
Here we need to note 
 that the cut cases
may cause a disconnected covering of $\bP^1$. 
Recall that the \emph{exponential}
generating function
$$
e^{\mathbf{H}(s,\mathbf{p})}
= 1 + {\mathbf{H}(s,\mathbf{p})} + 
\frac{1}{2}\;{\mathbf{H}(s,\mathbf{p})}^2
+
\frac{1}{3!}\;{\mathbf{H}(s,\mathbf{p})}^3+\cdots
$$
counts disconnected Hurwitz coverings.
The power of ${\mathbf{H}(s,\mathbf{p})}$ is the number
of connected components.  
Now the above merging consideration gives the following equation,
which is 
 the cut-and-join equation as a linear partial differential
 equation
\begin{equation}
\label{eq:caj-disconnected}
\left[\frac{\partial}{\partial s}-\frac{1}{2}\sum_{\alpha,\beta\ge 1}
\left((\alpha+\beta)p_\alpha p_\beta \frac{\partial}{\partial p_{\alpha+\beta}} + 
\alpha \beta p_{\alpha+\beta}\frac{\partial^2}{\partial p_\alpha\partial p_\beta}\right)\right]
 e^{\mathbf{H}(s,\mathbf{p})} = 0 .
\end{equation}
We can 
immediately deduce
\begin{equation}
\label{eq:caj-connected}
\frac{\partial \mathbf{H}}{\partial s}
=\frac{1}{2}\sum_{\alpha, \beta\ge 1}
\left((\alpha+\beta)p_\alpha p_\beta \frac{\partial \mathbf{H}}{\partial p_{\alpha+\beta}} +
\alpha\beta p_{\alpha+\beta}\frac{\partial^2 \mathbf{H}}{\partial p_\alpha\partial p_\beta}
+ 
\alpha\beta p_{\alpha+\beta}\frac{\partial \mathbf{H}}{\partial p_\alpha}\cdot 
\frac{\partial \mathbf{H}}{\partial p_\beta}\right) .
\end{equation}
This is the cut-and-join equation for the generating function
$\mathbf{H}(s,\mathbf{p})$ of the number of 
\emph{connected} Hurwitz coverings.

At this stage, we apply the ELSV formula (\ref{eq:ELSV}) to 
(\ref{eq:anotherHurwitz}). 
For a partition $\mu$ of length $\ell$, we define
\begin{equation}
\label{eq:Hg(mu)}
H_g(\mu)= \frac{|\Aut(\mu)|}{r(g,\mu)!}\cdot h_{g,\mu}
=
\sum_{n_1+\cdots+n_\ell \le 3g-3+\ell} \la \tau_{n_1}\cdots \tau_{n_\ell}\Lambda_g^{\vee}(1)\ra\,\,\prod_{i=1}^{\ell}\frac{\mu_i^{\mu_i+n_i}}{\mu_i!}.
\end{equation}
Then we have
\begin{equation}
\label{eq:Hgell-bold}
\mathbf{H}_{g,\ell}(s,\mathbf{p}) 
=\sum_{\mu:\ell(\mu)=\ell}\frac{1}{|\Aut(\mu)|}H_g(\mu)
\mathbf{p}_\mu s^{r(g,\mu)}
= \frac{1}{\ell!}\sum_{(\mu_1,\dots,\mu_\ell)\in \bN^\ell} 
H_g(\mu) \mathbf{p}_\mu s^{r(g,\mu)}  .
\end{equation}
The automorphism factor $|\Aut(\mu)|$
in the formula comes from the
re-summation.  For any function $f(\mu)$ in $\mu$, we
have a change of summation formula
\begin{equation}
\label{eq:autofactor}
\sum_{\mu\in\bN^\ell}f(\mu)
=
\sum_{\mu:\ell(\mu) = \ell}\frac{1}{\big|\Aut(\mu)\big|}
\sum_{\sigma\in S_\ell}f(\mu^\sigma),
\end{equation}
where $S_\ell$ is the permutation group of $\ell$ letters
and 
$$
\mu^\sigma = \left(\mu_{\sigma(1)},\dots,\mu_{\sigma(\ell)}
\right)\in \bN^\ell
$$ 
is the integer vector obtained by
permuting the parts of $\mu$ by $\sigma\in S_\ell$. 
If  $f(\mu)$ is a symmetric function, then the summation over
$S_\ell$ simply contributes $\ell!$ to the formula, 
as in (\ref{eq:Hgell-bold}).
For a partition $\mu$, let us denote by 
 $m_{\a}(\mu)$ the multiplicity of $\a$ in $\mu$, 
 i.e., the number of $\a$ repeated in $\mu$.
Then we have
\begin{equation}
\label{eq:Aut-mu}
\big|\Aut(\mu)\big| = \prod_{k\ge 1} m_k (\mu)!.
\end{equation}

Let us now compare the coefficient of 
$\mathbf{p}_\mu s^{r-1}$ in the cut-and-join equation
(\ref{eq:caj-connected}) for a given partition $\mu$ and an integer
$r\ge 1$.
The left-hand side contributes 
\begin{equation}
\label{eq:caj-LHS}
r(g,\mu)\;\frac{{H}_g  (\mu) }{|\Aut(\mu)|} ,
\end{equation}
subject to the condition $r=r(g,\mu)$.

The terms of $\mathbf{p}_\mu s^{r-1}$ that come from the
 \emph{cut}-operation of the right-hand side of (\ref{eq:caj-connected}) 
must have a profile $\big(\mu(\hat{i},\hat{j}),\mu_i+\mu_j\big)$, because 
\begin{align*}
r\left(g,\big(\mu(\hat{i},\hat{j}),\mu_i+\mu_j\big)\right)
&=
2g-2+\ell\big(\mu(\hat{i},\hat{j}),\mu_i+\mu_j\big)+\left|\big(\mu(\hat{i},\hat{j}),\mu_i+\mu_j\big)\right|
\\
&=
2g-2+(\ell-1)+|\mu|
 = r(g,\mu)-1.
\end{align*}
We see that the application of the differential operator
$p_{\mu_i}p_{\mu_j}\partial\big/\partial p_{\mu_i+\mu_j}$
to $\mathbf{H}(s,\mathbf{p})$
restores the profile $\mu$ from $\big(\mu(\hat{i},\hat{j}),\mu_i+\mu_j\big)$. 
Thus the coefficient of $\mathbf{p}_\mu s^{r-1}$ is 
\begin{equation}
\label{eq:cut}
\frac{1}{\big|\Aut(\mu)\big|}
 \sum_{i< j}
( \mu_i+\mu_j)
{H}_g \big(\mu(\hat{i},\hat{j}),\mu_i+\mu_j\big)
.
 \end{equation}
 In this consideration, we are naming all parts of
 $\mu$ to apply the cut-operation. Therefore, we need to compensate
 the overcount by the $\Aut(\mu)$-factor. 
 In terms of combinatorics, we can obtain (\ref{eq:cut})
 in a different way. It is easy to see  
 \cite[Section~2.3]{Zhou} that
 \begin{equation}
 \label{eq:Aut-cut}
 \left|
 \Aut\big(\mu(\hat{i},\hat{j}),\mu_i+\mu_j\big)
 \right|
 =
 \begin{cases}
  \big|
 \Aut(\mu)
 \big|
 \cdot
 \frac{m_{\mu_i+\mu_j}(\mu)+1}{m_{\mu_i}(\mu)m_{\mu_j}(\mu)} \qquad \mu_i\ne \mu_j,
 \\
  \big|
 \Aut(\mu)
 \big|
 \cdot
 \frac{m_{\mu_i+\mu_j}(\mu)+1}{m_{\mu_i}(\mu)\big(m_{\mu_i}(\mu)-1\big)} \qquad \mu_i= \mu_j.
 \end{cases}
 \end{equation}
So if $\mu_i\ne \mu_j$, then 
\begin{multline*}
\frac{1}{\big|\Aut(\mu)\big|}
( \mu_i+\mu_j)
{H}_g \big(\mu(\hat{i},\hat{j}),\mu_i+\mu_j\big)
\\
=
( \mu_i+\mu_j)
 \cdot
 \frac{m_{\mu_i+\mu_j}(\mu)+1}{m_{\mu_i}(\mu)m_{\mu_j}(\mu)}
 \cdot
\frac{
{H}_g \big(\mu(\hat{i},\hat{j}),\mu_i+\mu_j\big)
}
{ \left|
 \Aut\big(\mu(\hat{i},\hat{j}),\mu_i+\mu_j\big)
 \right|},
\end{multline*}
where each factor of the right-hand side has combinatorial significance.
When $\mu_i=\mu_j=\a$, we have
\begin{equation*}
\frac{1}{\big|\Aut(\mu)\big|}
( \mu_i+\mu_j)
{H}_g \big(\mu(\hat{i},\hat{j}),\mu_i+\mu_j\big)
=
\a
 \cdot
 \frac{m_{2\a}(\mu)+1}{
 \binom{m_{\a}(\mu)}{2}}
 \cdot
\frac{
{H}_g \big(\mu(\hat{i},\hat{j}),2\a\big)
}
{ \left|
 \Aut\big(\mu(\hat{i},\hat{j}),2\a\big)
 \right|},
\end{equation*}
where the part $\a$ is removed from the $i$-th and $j$-th slots
of $\mu$.

 In a \emph{join} term we must have a profile
$\big(\mu(\hat{i}),\a,\b\big)$. Since $\ell\big(\mu(\hat{i}),\a,\b\big) = \ell+1$,
changing $r$ to $r-1$ requires reducing the genus. 
One possibility is 
\begin{align*}
r\left(g-1,\big(\mu(\hat{i}),\a,\b\big)\right)
&=
2(g-1) -2 +\ell\big(\mu(\hat{i}),\a,\b\big) 
+\left|\big(\mu(\hat{i}),\a,\b\big)\right|
\\
&=
2g - 2 +(\ell + 1) +|\mu| -2= r(g,\mu)-1.
\end{align*}
In this case the differential operator 
$p_{\a +\b} \partial^2\big/\partial p_\a \partial p_\b$
applied to $\mathbf{H}(s,\mathbf{p})$ recovers the
profile $\mu$. The coefficient of $\mathbf{p}_\mu s^{r-1}$ is 
then
\begin{equation}
\label{eq:join1}
\frac{1}{2\big|\Aut(\mu)\big|} \sum_{i=1} ^\ell
 \sum_{\alpha+\beta = \mu_i}\alpha\beta 
{H}_{g-1} \big(\mu(\hat{i}),\a,\b\big)
.
\end{equation}
Here again we can give a combinatorial explanation
of this formula
using (\ref{eq:Aut-mu}) and (\ref{eq:Aut-cut}). 
When $\a\ne \b$, we have
\begin{equation*}
\frac{1}{\big|\Aut(\mu)\big|}
\alpha\beta 
{H}_{g-1} \big(\mu(\hat{i}),\a,\b\big)
=
\alpha\beta \cdot
 \frac{\big(m_\a(\mu)+1\big)\big(m_\b(\mu)+1\big)}
 {m_{\mu_i}(\mu)}
 \cdot
\frac{{H}_{g-1} \big(\mu(\hat{i}),\a,\b\big)}
{\left|
\Aut \big(\mu(\hat{i}),\a,\b\big)
\right|}.
\end{equation*}
And if $\a=\b=\half \mu_i$, then 
\begin{equation*}
\frac{1}{\big|\Aut(\mu)\big|} 
\alpha^2
{H}_{g-1} \big(\mu(\hat{i}),\a,\a\big)
=
2\alpha^2 \cdot
 \frac{\binom{m_\a(\mu)+2}{2}}
 {m_{\mu_i}(\mu)}
 \cdot
\frac{{H}_{g-1} \big(\mu(\hat{i}),\a,\a\big)}
{\left|
\Aut \big(\mu(\hat{i}),\a,\a\big)
\right|}.
\end{equation*}
The overall factor $2$ in the right-hand side comes from the
second order
differentiation $\partial^2\big/\partial p_\a ^2$.

There is yet another possibility
to obtain the profile $\mu$ from a \emph{join}-operation,
  if we utilize
disconnected Hurwitz covers.
Consider Hurwitz covers
$$
f_1:X_1\longrightarrow \bP^1
\qquad \text{and}\qquad 
f_2:X_2\longrightarrow \bP^1
$$
of genus $g_1$ (\emph{resp.} $g_2$) and profile $(\nu_1,\a)$ 
(\emph{resp.} $(\nu_2,\b)$).
Let  $\nu_1\sqcup\nu_2$ denote the partition obtained
by gathering all parts 
of $\nu_1$ and $\nu_2$ together.
If $g_1+g_2=g$ and  
$
\nu_1\sqcup \nu_2 = \mu(\hat{i})
$,
 then the join-operation 
recovers the profile $\mu$, provided that $\a+\b=\mu_i$. 
This is because
\begin{align*}
r\big(g_1,(\nu_1,\a)\big) &=2g_1-2+\ell(\nu_1)+1+|\nu_1|+\a
\\
r\big(g_2,(\nu_2,\b)\big)&=2g_2-2+\ell(\nu_2)+1+|\nu_2|+\b
\\
r(g,\mu)-1&=2g-2+\ell+|\mu|-1.
\end{align*}
The $\mathbf{p}_\mu s^{r-1}$-term comes from
$p_{\alpha+\beta}\frac{\partial \mathbf{H}}{\partial p_\alpha}\cdot 
\frac{\partial \mathbf{H}}{\partial p_\beta}$,
and its coefficient is 
\begin{equation}
\label{eq:join2}
\frac{1}{2\big|\Aut(\mu)\big|} \sum_{i=1} ^\ell
 \sum_{\alpha+\beta = \mu_i}\alpha\beta 
\sum_{\substack{g_1+g_2 = g\\
\nu_1\sqcup \nu_2 = \mu(\hat{i})}}
{H}_{g_1} (\nu_1,\a)
{H}_{g_2} (\nu_2,\b).
\end{equation}
The combinatorial derivation of this formula
follows from the identity
\begin{equation*}
\left|
\Aut(\nu_1\sqcup\nu_2)
\right|
=
\big|
\Aut(\nu_1)
\big|
\cdot
\big|
\Aut(\nu_2)
\big|
\cdot
\prod_{k\ge 1}\binom{m_k(\nu_1\sqcup\nu_2)}{m_k(\nu_1)}.
\end{equation*}
When $\a\ne \b$, we have
\begin{multline*}
\frac{1}{\big|\Aut(\mu)\big|} 
\alpha\beta 
\sum_{\substack{g_1+g_2 = g\\
\nu_1\sqcup \nu_2 = \mu(\hat{i})}}
{H}_{g_1} (\nu_1,\a)
{H}_{g_2} (\nu_2,\b)
\\
=
\alpha\beta \cdot
 \frac{\big(m_\a(\mu)+1\big)\big(m_\b(\mu)+1\big)}
 {m_{\mu_i}(\mu)}
 \cdot
 \frac{1}{\prod_{k\ge 1}\binom{m_k(\mu(\hat{i}),\a,\b)}{m_k(\nu_1,\a)}}
 \cdot
 \frac{{H}_{g_1} (\nu_1,\a)}{\big|\Aut(\nu_1,\a)\big|}
 \cdot
\frac{{H}_{g_2} (\nu_2,\b)}{\big|\Aut(\nu_2,\b)\big|}.
\end{multline*}
And if  $\a= \b=\half \mu_i$, then
\begin{multline*}
\frac{1}{\big|\Aut(\mu)\big|} 
\alpha^2
\sum_{\substack{g_1+g_2 = g\\
\nu_1\sqcup \nu_2 = \mu(\hat{i})}}
{H}_{g_1} (\nu_1,\a)
{H}_{g_2} (\nu_2,\a)
\\
=
2\alpha^2 \cdot
 \frac{\binom{m_\a(\mu)+2}{2}}
 {m_{\mu_i}(\mu)}
 \cdot
 \frac{1}{\prod_{k\ge 1}\binom{m_k(\mu(\hat{i}),\a,\a)}{m_k(\nu_1,\a)}}
 \cdot
 \frac{{H}_{g_1} (\nu_1,\a)}{\big|\Aut(\nu_1,\a)\big|}
 \cdot
\frac{{H}_{g_2} (\nu_2,\a)}{\big|\Aut(\nu_2,\a)\big|}.
\end{multline*}

Assembling (\ref{eq:cut}),
(\ref{eq:join1}), and (\ref{eq:join2})  together, we obtain
the combinatorial form of the cut-and-join equation.

\begin{thm}[Cut-and-join equation]
\label{thm:caj}
The functions $H_g(\mu)$ of {\rm{(\ref{eq:Hg(mu)})}}
satisfy a recursion equation
\begin{multline}
\label{eq:caj}
r(g,\mu){H}_g  (\mu)
=
\sum_{i< j}
(\mu_i+\mu_j)
{H}_g \big(\mu(\hat{i},\hat{j}),\mu_i+\mu_j\big)
 \\
 + \frac{1}{2} \sum_{i=1} ^\ell
 \sum_{\alpha+\beta = \mu_i}\alpha\beta 
\left({H}_{g-1} \big(\mu(\hat{i}),\a,\b\big)
+
\sum_{\substack{g_1+g_2 = g\\
\nu_1\sqcup \nu_2 = \mu(\hat{i})}}
{H}_{g_1} (\nu_1,\a)
{H}_{g_2} (\nu_2,\b)
\right).
\end{multline}
\end{thm}

\section{Laplace transform and the Lambert curve}
\label{sect:LT}

Since linear Hodge integrals
$\la \tau_{n_1}\cdots\tau_{n_\ell}\lam_j\ra$
do not depend on a partition $\mu$, 
it is natural to ask if there is any 
direct recursion formula for them \emph{without}
any reference to partitions. 
The answer is \emph{yes}, and we give the
formula in Section~\ref{sect:recursion}.
The natural complexity measure for the 
moduli space $\overline{\cM}_{g,\ell}$
is the absolute value $2g-2+\ell$ of the 
Euler characteristic of an $\ell$-punctured 
Riemann surface of genus $g$. 
An inductive formula associated to 
$\overline{\cM}_{g,\ell}$ with respect to 
$2g-2+\ell$ is generally called a \emph{topological recursion}.
We wish to establish a topological recursion 
for linear Hodge integrals. 
In the light of (\ref{eq:Hg(mu)}) and 
the combinatorial cut-and-join equation
(\ref{eq:caj}), it is obvious what we should do
to eliminate the $\mu$-dependence:
just take the summation over all partitions $\mu$. 
This is the idea of the \emph{Laplace transform}
discovered in \cite{EMS}. In this section we explain
this idea.

Since the sum of $\frac{k^{k+n}}{k!}$ for all positive
integer $k$ diverges, we are naturally led to the
idea of 
 Laplace transformation. Indeed,
\begin{equation}
\label{eq:fn}
f_n(w)=\sum_{k=1} ^\infty \frac{k^{k+n}}{k!}e^{-k(w+1)}
\end{equation}
is a holomorphic function in $w$ for $Re(w)>0$.
This follows from 
 Stirling's formula
$$
e^{-k}\frac{k^{k+n}}{k!}\sim \frac{1}{\sqrt{2\pi}}\;k^{n-\half}
\qquad {\text{for}}\qquad k>\!>1.
$$
Note that the continuous estimate for (\ref{eq:fn}) is given by
 \begin{equation}
 \label{eq:Stirling}
 \int_0 ^\infty x^{n-\half} e^{-xw} dx = \frac{\Gam(n+\half)}
 {w^{n+\half}}
 \end{equation}
for $n> -\half$. Thus $f_n(w)$ of (\ref{eq:fn}) is 
expected to be a function of $\sqrt{w}$, instead of 
$w$ itself,
if $n$ is an integer. We now come to the point of asking:
\emph{what is the Riemann surface of the function}
$f_n(w)$?
If the estimate (\ref{eq:Stirling}) were exact, then 
the Riemann surface of $f_n(w)$ would have been the same as 
that of $\sqrt{w}$. But since it is not, we need a
different idea.

The idea used in \cite{EMS} is the following. 
First we introduce a function
\begin{equation}
\label{eq:t}
t = t(w) = 1+ \sum_{k=1} ^\infty \frac{k^{k}}{k!}e^{-k(w+1)},
\end{equation}
which is holomorphic for $Re(w)>0$, and define
\begin{equation}
\label{eq:xy}
x = e^{-(w+1)}\qquad{\text{and}}\qquad
y=\frac{t-1}{t}.
\end{equation}
We can solve $t=t(w)$ in terms of $x$ and $y$.
The result is 
\begin{equation}
\label{eq:Lambert}
x = ye^{-y}.
\end{equation}
Let us call the plane analytic curve
\begin{equation}
\label{eq:Lambertcurve}
C = \{(x,y)\in\bC^2\;|\;x = ye^{-y}\}\subset \bC^2
\end{equation}
the \emph{Lambert curve}. This naming is
due to the resemblance of
(\ref{eq:t}) and the classical
\emph{Lambert W-function}
$$
W(x) = -\sum_{k=1}^{\infty}
\frac{k^{k-1}}{k!}(-x)^{k}.
$$
The Lambert curve $C$ is analytically isomorphic
to $\bC$, so it is an open Riemann surface of genus $0$. 
The $x$-projection $\pi:C\rightarrow \bC$
has a unique critical point $q_0=(e^{-1},1)\in C$. 
In terms of the coordinates $w$ and $t$, the inverse function
of (\ref{eq:t}), or the equation for the 
Lambert curve, is given by
\begin{equation}
\label{eq:w}
 w =w(t)= -\frac{1}{t} - \log \left(1-\frac{1}{t}\right)
=\sum_{m=2} ^\infty \frac{1}{m}\;\frac{1}{t^m} ,
\end{equation}
which is holomorphic for $Re(t)>1$. 
The critical point of the projection $\pi$ in this coordinate
is $(w,t) = (0,\infty)$. Since the infinite series of (\ref{eq:w})
starts at $m=2$, $\pi$ is locally a double-sheeted covering
around $w=0$. And this is what we wanted.
Indeed, the Lambert curve $C$ is \emph{the Riemann surface}
of the function $f_n(w)$. 
It is natural to consider $f_n(w)$ as a function
in $t$, since $t$ is a global coordinate of $C$.
So we re-define
\begin{equation}
\label{eq:hxin}
\hxi_n(t)=\sum_{k=1} ^\infty \frac{k^{k+n}}{k!}e^{-k(w+1)},
\end{equation}
which is simply $f_n(w)$ in terms of $t$ satisfying 
$w=w(t)$. 
But something remarkable happens here: $\hxi_n(t)$ is a 
\emph{polynomial} in $t$ if $n\ge 0$.
The proof is  obvious. 
A standard property of the Laplace transform gives
\begin{equation}
\label{eq:LTdiff}
-\frac{d}{dw}f_n(w) 
= \sum_{k=1} ^\infty \frac{k^{k+n+1}}{k!}e^{-k(w+1)}
=f_{n+1}(w),
\end{equation}
and the coordinate change (\ref{eq:w}) implies
\begin{equation}
\label{eq:w-in-t}
-\frac{d}{dw} = t^2(t-1) \frac{d}{dt}.
\end{equation}
Therefore, $\hxi_n(t)$'s satisfy a recursion formula
\begin{equation}
\label{eq:hxi-recursion}
\hxi_{n+1}(t) = t^2(t-1) \frac{d}{dt}\hxi_n(t)
=D\hxi_n(t).
\end{equation}
Since $\hxi_0(t) = t-1$ from (\ref{eq:t}), we see that 
$\hxi_n(t)$ is a polynomial in $t$ of degree $2n+1$. 
It immediately follows that the 
Laplace transform  
\begin{multline}
\label{eq:Hhatgell}
\hatH_{g,\ell}(t_1,\dots,t_\ell) = 
\sum_{\mu\in\bN^\ell} H_g(\mu)
e^{-\left(
\mu_1(w_1+1)+\cdots+\mu_\ell(w_\ell+1)
\right)}
\\
=
\sum_{n_1+\cdots+n_\ell \le 3g-3+\ell} \la \tau_{n_1}\cdots \tau_{n_\ell}\Lambda_g^{\vee}(1)\ra\,\,\prod_{i=1}^{\ell}
\hxi_{n_i}(t_i)
\end{multline}
of $H_g(\mu)$ of (\ref{eq:Hg(mu)}) 
is  a symmetric \emph{polynomial} in the $t$-variables
 and naturally lives 
on $C^\ell$, when $2g-2+\ell>0$.

The unstable geometries $(g,\ell) = (0,1)$ and 
$(0,2)$ are the exceptions for this general formula.
Recall the $(0,1)$ case (\ref{eq:01Hodge}). 
We have
\begin{equation}
\label{eq:H01LT}
\hatH_{0,1}(t) = \sum_{k=1} ^\infty \frac{k^{k-2}}{k!}
e^{-k(w+1)} = -\frac{1}{2 \,t^2}+c = \hxi_{-2}(t),
\end{equation}
where the constant $c$ is given  by
$$
c=\sum_{k=1} ^\infty \frac{k^{k-2}}{k!}e^{-k}.
$$
The $(0,2)$ case (\ref{eq:02Hodge}) is quite more involved.
It is proved  \cite[Proposition~3.6]{EMS}  that we have
\begin{multline}
\label{eq:H02LT}
\widehat{\cH}_{0,2}(t_1,t_2) =
\sum_{\mu_1,\mu_2\ge 1}
\frac{1}{\mu_1+\mu_2}\cdot\frac{\mu_1 ^{\mu_1}}{\mu_1 !}\cdot\frac{\mu_2 ^{\mu_2}}{\mu_2 !}
e^{-\mu_1 (w_1+1)}e^{-\mu_2 (w_2+1)}\\
=\log\left(
\frac{\hxi_{-1}(t_1)-\hxi_{-1}(t_2)}
{x_1-x_2}
\right)
-\hxi_{-1}(t_1)-\hxi_{-1}(t_2),
\end{multline}
where
\begin{equation}
\label{eq:hxi-1}
\hxi_{-1}(t) = \frac{t-1}{t} = y.
\end{equation}

\section{The topological recursion as a Laplace transform} \label{sect:recursion}

In the previous section we have 
computed the Laplace transform of $H_g(\mu)$ as
a function on partitions $\mu$.  In this section we
calculate the Laplace transform of the
cut-and-join equation (\ref{eq:caj}) and prove
Theorem~\ref{thm:main}.

Let us denote 
\begin{equation}
\label{eq:linearform}
\la \mu,w+1\ra 
= \mu_1(w_1+1)+\cdots+
\mu_\ell(w_\ell+1).
\end{equation}
Recalling the expression of $r(g,\mu)$ given in (\ref{eq:RH})
and using (\ref{eq:autofactor}), 
the Laplace transform of the left-hand side of (\ref{eq:caj}) becomes
\begin{equation}
\label{eq:LTcaj-LHS}
\sum_{\mu\in\bN^\ell}r(g,\mu)
H_g(\mu)e^{-\la\mu,w+1\ra}
\\
=
\left(
2g-2+\ell +\sum_{i=1} ^\ell
t_i ^2 (t_i-1)\frac{\partial }{\partial t_i}
\right)
\hatH_{g,\ell}(t_1,\dots,t_\ell).
\end{equation}
Here we note that multiplication of $\mu_i$ to the summand corresponds
to the operation of $D_i=t_i ^2(t_i-1)\frac{\partial }{\partial t_i}$
due to (\ref{eq:w-in-t}).

  To find  the Laplace transform of the cut terms (\ref{eq:cut}),
we first note a formula:
 \begin{align*}
 \sum_{\mu_1,\mu_2\ge 0}f(\mu_1+\mu_2)e^{-(\mu_1w_1+
 \mu_2 w_2)}
 &=
 \sum_{k=0} ^\infty \sum_{m=0} ^{k}
 f(k)e^{-kw_1}e^{-m(w_2-w_1)}
 \\
&=
 \sum_{k=0} ^\infty \frac{1-e^{-(k+1)(w_2-w_1)}}
 {1-e^{-(w_2-w_1)}}\;f(k)e^{-kw_1}
 \\
&=
 \frac{1}{e^{-w_1}-e^{-w_2}}
 \sum_{k=0} ^\infty f(k)\left(
 e^{-(k+1)w_1}-e^{-(k+1)w_2}
 \right).
 \end{align*}
 Thus we obtain
 \begin{align*}
 &\qquad
 \half 
 \sum_{\mu\in\bN^\ell}
 \sum_{i \ne j}
 (\mu_i+\mu_j)H_g\left(\mu(\hat{i},\hat{j}),\mu_i+\mu_j\right)
 e^{-\la\mu,w+1\ra}
 \\
 &=
 \half
\sum_{i \ne j}
 \frac{1}{e^{-(w_i+1)}-e^{-(w_j+1)}}
    \Bigg(
 e^{-(w_i+1)} 
   t_i ^2(t_i-1)\frac{\partial}{\partial t_i}
    \hatH_{g,\ell-1}\left(t_1,\dots,\widehat{t_j},\dots,t_\ell\right)
    \\
 &\qquad-
  e^{-(w_j+1)} 
   t_j ^2(t_j-1)\frac{\partial}{\partial t_j}
    \hatH_{g,\ell-1}\left(t_1,\dots,\widehat{t_i},\dots,t_\ell\right)
  \Bigg)
   \\
&\qquad- 
\sum_{i\ne j} t_i ^2(t_i-1)\frac{\partial}{\partial t_i}
\hatH_{g,\ell-1}\left(t_1,\dots,\widehat{t_j},\dots,t_\ell\right),
\end{align*}
where the last term comes from the adjustment of 
the cases $\mu_i=0$ and $\mu_j=0$ that are
not included in the Laplace transform.

The Laplace transform of the
first join terms (\ref{eq:join1})  is given by
\begin{multline*}
\half
\sum_{\mu\in\bN^\ell}
\sum_{i=1} ^\ell
\sum_{\a+\b=\mu_i}
\a\b {H}_{g-1} \big(\mu(\hat{i}),\a,\b\big)
e^{-\la \mu,w+1\ra}
\\
=
\sum_{i=1} ^\ell
\left[
u_1 ^2 (u_1-1) u_2 ^2(u_2-1)
\frac{\partial ^2}{\partial u_1\partial u_2}
\hatH_{g-1,\ell+1}\left(u_1,u_2,t_{L\setminus\{i\}}\right)
\right]_{u_1=u_2=t_i} ,
\end{multline*}
where  $t_I=(t_i)_{i\in I}$ for a subset 
$I\subset L=\{1,2,\dots,\ell\}$.
In the same way we can calculate the Laplace transform
of the second join terms (\ref{eq:join2}):
\begin{multline*}
\sum_{\mu\in\bN^\ell}
\sum_{\a+\b=\mu_i}
\a\b
\sum_{\substack{g_1+g_2 = g\\
\nu_1\sqcup \nu_2 = \mu(\hat{i})}}
{H}_{g_1} (\nu_1,\a)
{H}_{g_2} (\nu_2,\b)
e^{-\la \mu,w+1\ra}
\\
=
\left[
\sum_{\substack{g_1+g_2 = g\\
J\sqcup K= L\setminus\{i\}}}
u_1 ^2 (u_1-1) 
\frac{\partial }{\partial u_1}
\hatH_{g_1,|J|+1}(u_1,t_J)
u_2 ^2(u_2-1)
\frac{\partial}{\partial u_2}
\hatH_{g_2,|K|+1}(u_2,t_K)
\right]_{u_1=u_2=t_i}.
 \end{multline*}
Thus we establish
\begin{multline}
\label{eq:cajLT-pre}
\left(
2g-2+\ell +\sum_{i=1} ^\ell
t_i ^2 (t_i-1)\frac{\partial }{\partial t_i}
\right)
\hatH_{g,\ell}(t_1,\dots,t_\ell)
\\
=
\sum_{i<j}
 \frac{1}{e^{-(w_i+1)}-e^{-(w_j+1)}}
   \Bigg(
 e^{-(w_i+1)} 
   t_i ^2(t_i-1)\frac{\partial}{\partial t_i}
    \hatH_{g,\ell-1}\left(t_1,\dots,\widehat{t_j},\dots,t_\ell\right)
    \\
 -
  e^{-(w_j+1)} 
   t_j ^2(t_j-1)\frac{\partial}{\partial t_j}
    \hatH_{g,\ell-1}\left(t_1,\dots,\widehat{t_i},\dots,t_\ell\right)
  \Bigg)
   \\
- \sum_{i\ne j} t_i ^2(t_i-1)\frac{\partial}{\partial t_i}
\hatH_{g,\ell-1}\left(t_1,\dots,\widehat{t_j},\dots,t_\ell\right)
\\
+
\sum_{i=1} ^\ell
\left[
u_1 ^2 (u_1-1) u_2 ^2(u_2-1)
\frac{\partial ^2}{\partial u_1\partial u_2}
\hatH_{g-1,\ell+1}\left(u_1,u_2,t_{L\setminus\{i\}}\right)
\right]_{u_1=u_2=t_i}
\\
+
\half
\sum_{i=1} ^\ell
\sum_{\substack{g_1+g_2 = g\\
J\sqcup K= L\setminus\{i\}}}
t_i ^2 (t_i-1) 
\frac{\partial }{\partial t_i}
\hatH_{g_1,|J|+1}(t_i,t_J)\cdot 
t_i ^2 (t_i-1) 
\frac{\partial }{\partial t_i}
\hatH_{g_2,|K|+1}(t_i,t_K) .
\end{multline}
Note that  \emph{unstable geometries} are contained in 
the last summation. We use (\ref{eq:H01LT}) and 
(\ref{eq:H02LT}) to substitute the values in (\ref{eq:cajLT-pre}).
The result becomes surprisingly simple due to cancellation
of the non-polynomial terms. 
For $g_1=0$ and $J=\emptyset$, 
the contribution is
$$
\sum_{i=1} ^\ell
\hxi_{-1}(t_i)
t_i ^2 (t_i-1) 
\frac{\partial }{\partial t_i}
\hatH_{g,\ell}(t_1,\dots,t_\ell).
$$
For $g_1=0$ and $J=\{j\}\subset L\setminus\{i\}$, we have
\begin{multline*}
t_i ^2 (t_i-1) 
\frac{\partial }{\partial t_i}
\hatH_{0,2}(t_i,t_j)
=
\frac{\hxi_0(t_i)}{\hxi_{-1}(t_i)-\hxi_{-1}(t_j)}
-
\frac{x_i}{x_i-x_j}
-\hxi_0(t_i)
\\
=
\frac{\hxi_0(t_i)}{\hxi_{-1}(t_i)-\hxi_{-1}(t_j)}
-
\frac{e^{-(w_i+1)}}{e^{-(w_i+1)}-e^{-(w_j+1)}}
-\hxi_0(t_i).
\end{multline*}
Thus the unstable $(0,2)$ contribution in (\ref{eq:cajLT-pre})
is
\begin{multline*}
\sum_{i< j}
 \frac{  t_i ^2(t_i-1)^2\frac{\partial}{\partial t_i}
    \hatH_{g,\ell-1}\left(t_1,\dots,\widehat{t_j},\dots,t_\ell\right)
    -
      t_j ^2(t_j-1)^2\frac{\partial}{\partial t_j}
    \hatH_{g,\ell-1}\left(t_1,\dots,\widehat{t_i},\dots,t_\ell\right)}{\hxi_{-1}(t_i)-\hxi_{-1}(t_j)}
    \\
-
\sum_{i<j}
 \frac{1}{e^{-(w_i+1)}-e^{-(w_j+1)}}
   \Bigg(
 e^{-(w_i+1)} 
   t_i ^2(t_i-1)\frac{\partial}{\partial t_i}
    \hatH_{g,\ell-1}\left(t_1,\dots,\widehat{t_j},\dots,t_\ell\right)
    \\
 -
  e^{-(w_j+1)} 
   t_j ^2(t_j-1)\frac{\partial}{\partial t_j}
    \hatH_{g,\ell-1}\left(t_1,\dots,\widehat{t_i},\dots,t_\ell\right)
  \Bigg)
\\
-
\sum_{i\ne j}
\hxi_0(t_i)
 t_i ^2(t_i-1)\frac{\partial}{\partial t_i}
    \hatH_{g,\ell-1}\left(t_1,\dots,\widehat{t_j},\dots,t_\ell\right).
\end{multline*}
We have thus proved the following, which is 
equivalent to Theorem~\ref{thm:main}.

\begin{thm}
\label{thm:CAJLT}
The Laplace transform of the cut-and-join equation
is the following equation
for polynomials $\hatH_{g,\ell}(t_1,\dots,t_\ell)$
subject to the stability condition
$2g-2+\ell>0:$
\begin{multline}
\label{eq:cajLT}
\left(
2g-2+\ell +\sum_{i=1} ^\ell
\big(1-\hxi_{-1}(t_i)\big)
t_i ^2 (t_i-1)\frac{\partial }{\partial t_i}
\right)
\hatH_{g,\ell}(t_1,\dots,t_\ell)
\\
=
\sum_{i< j}t_it_j
 \frac{  t_i ^2(t_i-1)^2\frac{\partial}{\partial t_i}
    \hatH_{g,\ell-1}\left(t_1,\dots,\widehat{t_j},\dots,t_\ell\right)
    -
      t_j ^2(t_j-1)^2\frac{\partial}{\partial t_j}
    \hatH_{g,\ell-1}\left(t_1,\dots,\widehat{t_i},\dots,t_\ell\right)}{t_i-t_j}
   \\
- \sum_{i\ne j} t_i ^3(t_i-1)\frac{\partial}{\partial t_i}
\hatH_{g,\ell-1}\left(t_1,\dots,\widehat{t_j},\dots,t_\ell\right)
\\
+
\half
\sum_{i=1} ^\ell
\left[
u_1 ^2 (u_1-1) u_2 ^2(u_2-1)
\frac{\partial ^2}{\partial u_1\partial u_2}
\hatH_{g-1,\ell+1}\left(u_1,u_2,t_{L\setminus\{i\}}\right)
\right]_{u_1=u_2=t_i}
\\
+
\half
\sum_{i=1} ^\ell
\sum_{\substack{g_1+g_2 = g\\
J\sqcup K= L\setminus\{i\}}} ^{\rm{stable}}
t_i ^2 (t_i-1) 
\frac{\partial }{\partial t_i}
\hatH_{g_1,|J|+1}(t_i,t_J)\cdot 
t_i ^2 (t_i-1) 
\frac{\partial }{\partial t_i}
\hatH_{g_2,|K|+1}(t_i,t_K) .
\end{multline}
In the last sum each term is restricted to 
satisfying the stability conditions
$2g_1-1+|J|>0$ and $2g_2-1+|K|>0$.
\end{thm}

\begin{rem}
Eqn.(\ref{eq:cajLT}) is equivalent to the cut-and-join
equation (\ref{eq:caj-connected})
and (\ref{eq:caj}). Many other equivalent formulations
have been established, including the differential equation 
 of
\cite[Theorem~3.1]{GJV3}.
\end{rem}

\section{The Witten-Kontsevich
theorem and the $\lam_g$ formula}
\label{sect:WK}

It has been noticed that the asymptotic behavior of
Hurwitz numbers for a 
large partition recovers the intersection numbers of
$\psi$-classes \cite{OP1}. Actual recovery of the Witten-Kontsevich 
theorem \cite{K1992, W1991} from the ELSV formula
using this asymptotic argument is rather involved
(\cite{OP1}, see also \cite{KL}). 
Since the Laplace transform contains all the information 
of the asymptotics, we can easily deduce the
Virasoro constraint equation, or the equivalent Dijkgraaf-Verlinde-Verlinde
formula \cite[Eqn.~4.1]{DVV},
for the $\psi$-class intersection
 from our main equation (\ref{eq:main}). 
 Thus we obtain a straightforward proof of the Witten conjecture.
 In this section we observe that the \emph{top} degree terms
 of the recursion is the DVV formula. We also examine that
 the \emph{lowerst} degree terms imply the descendant relation
 of  the $\lam_g$ formula \cite{FP1, FP2}. 
 Our argument is along the same line with \cite{CLL, GJV3, Kazarian}. However,
 due to the polynomial formulation of (\ref{eq:main}), the
 derivation becomes simpler.

First we compute the polynomial $\hxi_n(t)$
using  (\ref{eq:hxi-recursion}). It has the general form
\begin{equation}
\label{eq:hxi-expansion}
\hxi_n(t)=
(2n-1)!! t^{2n+1} - \frac{(2n+1)!!}{3}\; t^{2n}+\cdots
+a_n t^{n+2} + (-1)^n n! \;t^{n+1},
\end{equation}
where $a_n$ is defined by
$$
a_{n} =-\big[(n+1)a_{n-1} +(-1)^n n!\big]
$$
and is identified as the sequence A001705 or A081047 of the
\emph{On-Line Encyclopedia of Integer Sequences}.

The DVV formula for the
Virasoro constraint condition on the $\psi$-class
intersections  reads 
\begin{multline}
\label{eq:DVV}
\la \tau_{n_L}\ra _{g,\ell}=
\sum_{j\ge 2} 
\frac{(2n_1+2n_j-1)!!}{(2n_1+1)!! (2n_j -1)!!}
 \la \tau_{n_1+n_j-1}
 \tau_{n_{L\setminus \{1,j\}}}
\ra _{g,\ell-1}
\\
+
\frac{1}{2} \sum_{a+b=n_1-2}
\left(
\la \tau_a\tau_b\tau_{n_{L\setminus\{1\}}}\ra _{g-1,\ell+1}
+
\sum_{\substack{g_1+g_2=g\\
J\sqcup K= L\setminus\{1\}}} ^{
\text{stable}}
\la \tau_a\tau_{n_J}\ra _{g_1,|J|+1}\cdot
\la \tau_b\tau_{n_K}\ra _{g_2,|K|+1}
\right)
\\
\times
\frac{(2a+1)!!(2b+1)!!}{(2n_1+1)!!}.
\end{multline}
Here $L=\{1,\dots,\ell\}$ is the index set as before, and
for a subset $I\subset L$ we write
$$
n_I = (n_i)_{i\in I}\qquad{\text{and}}\qquad
\tau_{n_I}=\prod_{i\in I}\tau_{n_i}.
$$

\begin{prop}
\label{prop:WK}
The DVV formula {\rm(\ref{eq:DVV})} is exactly the 
relation among the top degree coefficients
of the recursion {\rm(\ref{eq:main})}.
\end{prop}

\begin{proof}
Choose $n_L$ so that $|n_L|= n_1+n_2+\cdots +n_\ell=3g-3+\ell$.
The degree of the left-hand side of (\ref{eq:main}) is $3(2g-2+\ell) +1$. 
So we
 compare the coefficients of $t_1 ^{2n_1+2}\prod_{j\ge 2} t_j ^{2n_j+1}$
in the recursion formula. 
The contribution from the left-hand side of (\ref{eq:main}) is 
$$
\la \tau_{n_L}\ra _{g,\ell}
(2n_1+1)!!\prod_{j\ge 2} (2n_j-1)!!.
$$
The contribution from the first line of the right-hand side comes from
\begin{multline*}
\sum_{j\ge 2}
\la\tau_m\tau_{n_{L\setminus\{1,j\}}}\ra_{g,\ell-1}
(2m+1)!! 
 \frac{t_1^2 t_j t_1 ^{2m+3}
-t_j^2t_1t_j ^{2m+3}}{t_1-t_j}
\\
=
\sum_{j\ge 2}
\la\tau_m\tau_{n_{L\setminus\{1,j\}}}\ra_{g,\ell-1}
(2m+1)!! t_1t_j
 \frac{ t_1 ^{2m+4}
-t_j ^{2m+4}}{t_1-t_j}
\\
=
\sum_{j\ge 2}
\la\tau_m\tau_{n_{L\setminus\{1,j\}}}\ra_{g,\ell-1}
(2m+1)!!
\sum_{a+b=2m+3}t_1 ^{a+1}t_j ^{b+1},
\end{multline*}
where $m=n_1+n_j-1$. The matching term in this
formula is $a=2n_1+1$ and $b=2n_j$. Thus we extract
as the coefficient of $t_1 ^{2n_1+2}\prod_{j\ge 2} t_j ^{2n_j+1}$
$$
\sum_{j\ge 2}
\la\tau_{n_1+n_j-1}\tau_{n_{L\setminus\{1,j\}}}\ra_{g,\ell-1}
(2n_1+2n_j-1)!!\prod_{k\ne 1,j}(2n_k-1)!!.
$$
The contributions of the second and the third lines 
of the right-hand side  of (\ref{eq:main}) are
\begin{multline*}
\half
\sum_{a+b=n_1-2}
\left(
\la \tau_a \tau_b \tau_{L\setminus \{1\}}\ra_{g-1,\ell+1}
+\frac{1}{2}
\sum_{\substack{g_1+g_2=g\\
J\sqcup K= L\setminus\{1\}}} ^{
\text{stable}}
\la \tau_a\tau_{n_J}\ra _{g_1,|J|+1}\cdot
\la \tau_b\tau_{n_K}\ra _{g_2,|K|+1}
\right)
\\
\times
(2a+1)!!(2b+1)!!\prod_{j\ge 2}(2n_j-1)!!.
\end{multline*}
We have thus recovered the Witten-Kontsevich theorem
\cite{DVV, K1992, W1991}.
\end{proof}

The $\lam_g$ formula \cite{FP1, FP2, L, LLZ} is
\begin{equation}
\label{eq:lam-g}
\la \tau_{n_L}\lam_g\ra_{g,\ell}
=\binom{2g-3+\ell}{n_L}b_g,
\end{equation}
where 
\begin{equation}
\label{eq:multinomial}
\binom{2g-3+\ell}{n_L}=
\binom{2g-3+\ell}{n_1,\dots,n_\ell}
\end{equation}
is the multinomial coefficient, and 
$$
b_g = \frac{2^{2g-1}-1}{2^{2g-1}}\;
\frac{|B_{2g}|}{(2g)!}
$$
is a coefficient of the series
$$
\sum_{j=0} ^\infty b_j s^{2j}
= \frac{s/2}{\sin(s/2)}.
$$

\begin{prop}
\label{prop:lam-g}
The lowest degree terms of the topological recursion {\rm{(\ref{eq:main})}}
proves the combinatorial factor of the
$\lam_g$ formula
\begin{equation}
\label{eq:lam-g-combinatorial}
\la \tau_{n_L}\lam_g\ra_{g,\ell}
=\binom{2g-3+\ell}{n_L}
\la \tau_{2g-1}\lam_g\ra_{g,1}.
\end{equation}

\end{prop}

\begin{proof}
Choose  $n_L$ subject to $|n_L|=2g-3+\ell$. 
We compare the coefficient of  the terms of
 $\prod_{i\ge 1}
t_i ^{n_i+1}$ in (\ref{eq:main}), which has
degree $|n_L|+\ell =2g-3+2\ell$. The left-hand side contributes
\begin{multline*}
(-1)^{2g-3+\ell}(-1)^g
\la\tau_{n_L}\lam_g\ra_{g,\ell}
\prod_{i\ge 1} n_i!
\left(
2g-2+\ell - \sum_{i=1} ^\ell (n_i+1)
\right)
\\
=
(-1)^{\ell}(-1)^g
\la\tau_{n_L}\lam_g\ra_{g,\ell}
(\ell-1)
\prod_{i\ge 1}n_i!.
\end{multline*}
The lowest degree terms of the
first line of the right-hand side are
\begin{equation*}
(-1)^g
\sum_{i< j}\sum_m
\la \tau_m\tau_{L\setminus\{i,j\}}\lam_g\ra _{g,\ell-1}
(-1)^{m}(m+1) !
 \frac{  
 t_i ^{m+4}
    -
    t_j^{m+4}}{t_i-t_j}
    (-1)^{2g-3+\ell-n_i-n_j}\prod_{k\ne i,j}n_k! t_k ^{n_k+1}.
\end{equation*}
Since $m=n_i+n_j-1$, the coefficient of  $\prod_{i\ge 1}
t_i ^{n_i+1}$ is
\begin{equation*}
-(-1)^g(-1)^{2g-3+\ell}\sum_{i< j}
\la \tau_{n_i+n_j-1}\tau_{L\setminus\{i,j\}}\lam_g\ra _{g,\ell-1}
\binom{n_i+n_j}{n_i}
\prod_{i\ge 1}n_i! .
\end{equation*}
Note that the lowest degree coming from the second and the third
lines of the right-hand side of (\ref{eq:main}) is $|n_L|+\ell + 2$, which
is higher than the lowest degree of the left-hand side. Therefore, we have 
obtained a recursion equation with respect to $\ell$
\begin{equation}
\label{eq:lam-g-recursion}
(\ell-1)\la \tau_{n_L}\lam_g \ra_{g,\ell}=
\sum_{i< j}
\la \tau_{n_i+n_j-1}\tau_{L\setminus\{i,j\}}\lam_g\ra _{g,\ell-1}
\binom{n_i+n_j}{n_i}.
\end{equation}
The solution of the recursion equation (\ref{eq:lam-g-recursion})
is the multinomial coefficient (\ref{eq:multinomial}).
\end{proof}

\begin{rem}
Although the topological recursion 
(\ref{eq:main}) determines all linear Hodge integrals,
the closed formula
$$
b_g = \la \tau_{2g-2}\lam_g\ra_{g,1}\qquad g\ge 1
$$ 
does not seem to follow directly from it.
\end{rem}


\providecommand{\bysame}{\leavevmode\hbox to3em{\hrulefill}\thinspace}

\bibliographystyle{amsplain}

\end{document}